\documentclass[12pt]{article}
\usepackage[final]{epsfig}
\usepackage{graphics}
\usepackage{amsmath}
\usepackage{amsfonts}
\usepackage{latexsym}
\usepackage{amssymb}
\usepackage{graphicx}
\usepackage{url}
\usepackage{epstopdf}

\newtheorem{theorem}{Theorem}

\begin{document}
\newcommand{\eps}{{\varepsilon}}
\newcommand{\proofend}{$\Box$\bigskip}
\newcommand{\C}{{\mathbb C}}
\newcommand{\Q}{{\mathbb Q}}
\newcommand{\R}{{\mathbb R}}
\newcommand{\Z}{{\mathbb Z}}
\newcommand{\RP}{{\mathbb {RP}}}
\newcommand{\CP}{{\mathbb {CP}}}
\newcommand{\Tr}{\rm Tr}
\def\proof{\paragraph{Proof.}}

\title{Kasner meets Poncelet}

\author{Serge Tabachnikov\footnote{
Department of Mathematics,
Penn State University,
University Park, PA 16802;
tabachni@math.psu.edu}
}

\date{}
\maketitle

Given a planar pentagon $P$, construct two new pentagons, $D(P)$ and $I(P)$: the vertices of $D(P)$ are the intersection points of the diagonals of $P$, and the vertices of $I(P)$ are the tangency points of the conic inscribed in $P$ (the vertices and the tangency points are taken in their cyclic order). 
The following result is due to E. Kasner  \cite{Ka} (published in 1928, but discovered much earlier, in 1896).

\begin{figure}[hbtp]
\centering
\includegraphics[height=2in]{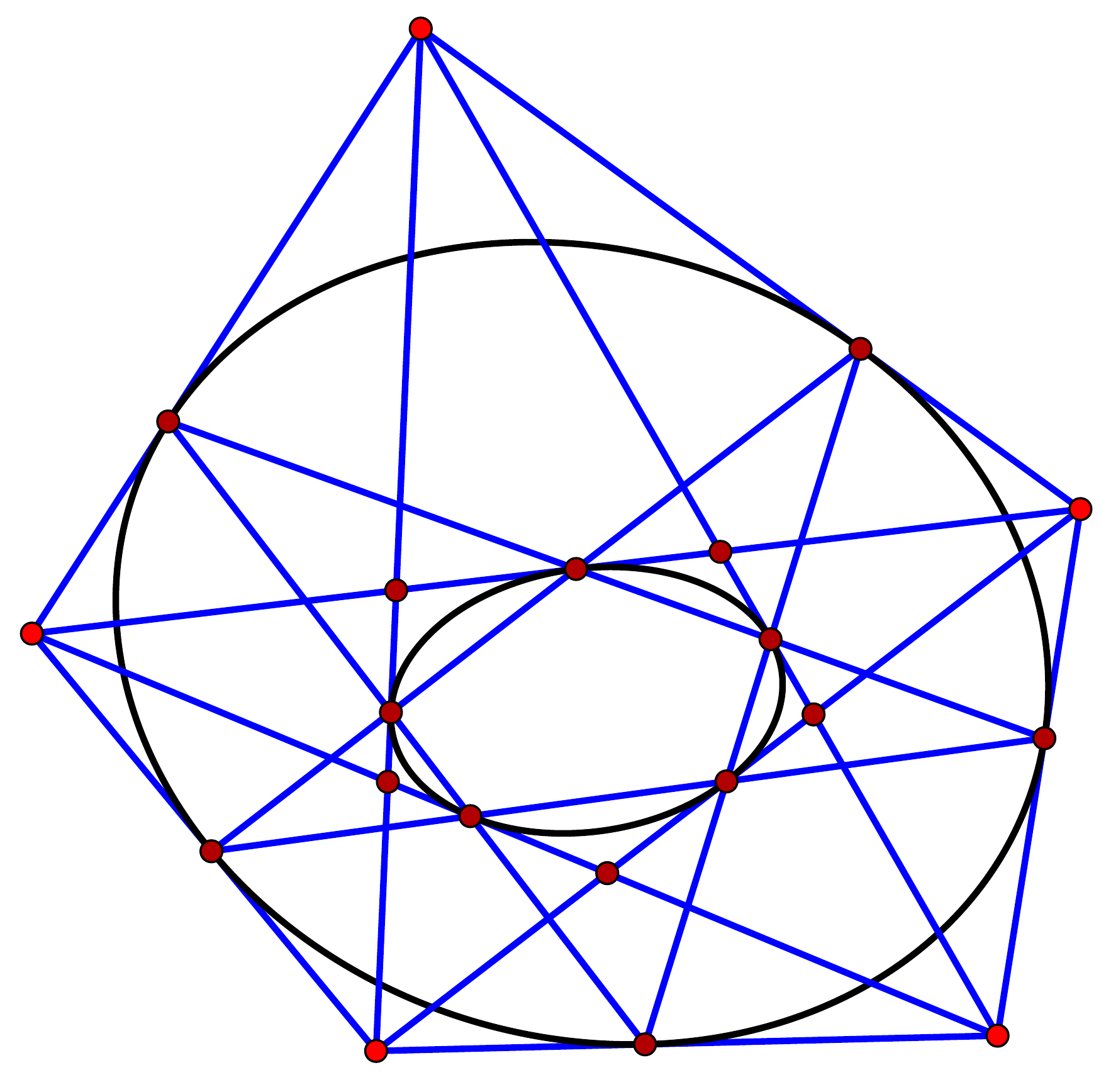}
\caption{Kasner's theorem.}
\label{penta}
\end{figure}

\begin{theorem} \label{Kas}
The two operations on pentagons, $D$ and $I$, commute: $ID(P)=DI(P)$, see Figure \ref{penta}.
\end{theorem}

A polygon inscribed into a conic and circumscribed about a conic is called a Poncelet polygon. The celebrated Poncelet porism states that the existence of a Poncelet $n$-gon on a pair of nested ellipses implies that every point of the outer ellipse is a vertex of such a Poncelet $n$-gon, see, e.g., \cite{DR,Fl}.

Since a (generic) quintuple of points lies on a conic, a (generic) quintuple of lines is tangent to a conic, a (generic) pentagon is a Poncelet polygon. We generalize Kasner's theorem  from pentagons to Poncelet polygons. 

Let $P$ be a Poncelet $n$-gon, $n\ge 5$. As before, the $n$-gon $I(P)$ is formed by the consecutive tangency points of the sides of $P$ with the inscribed conic. Fix $2\leq k < n/2$, and draw the $k$-diagonals of $P$, that is, connect $i$th vertex with $(i+k)$th vertex for $i=1,\ldots,n$. The consecutive intersection 
points of these diagonals form a new $n$-gon, denoted by $D_k (P)$. 

\begin{figure}[hbtp]
\centering
\includegraphics[height=4in]{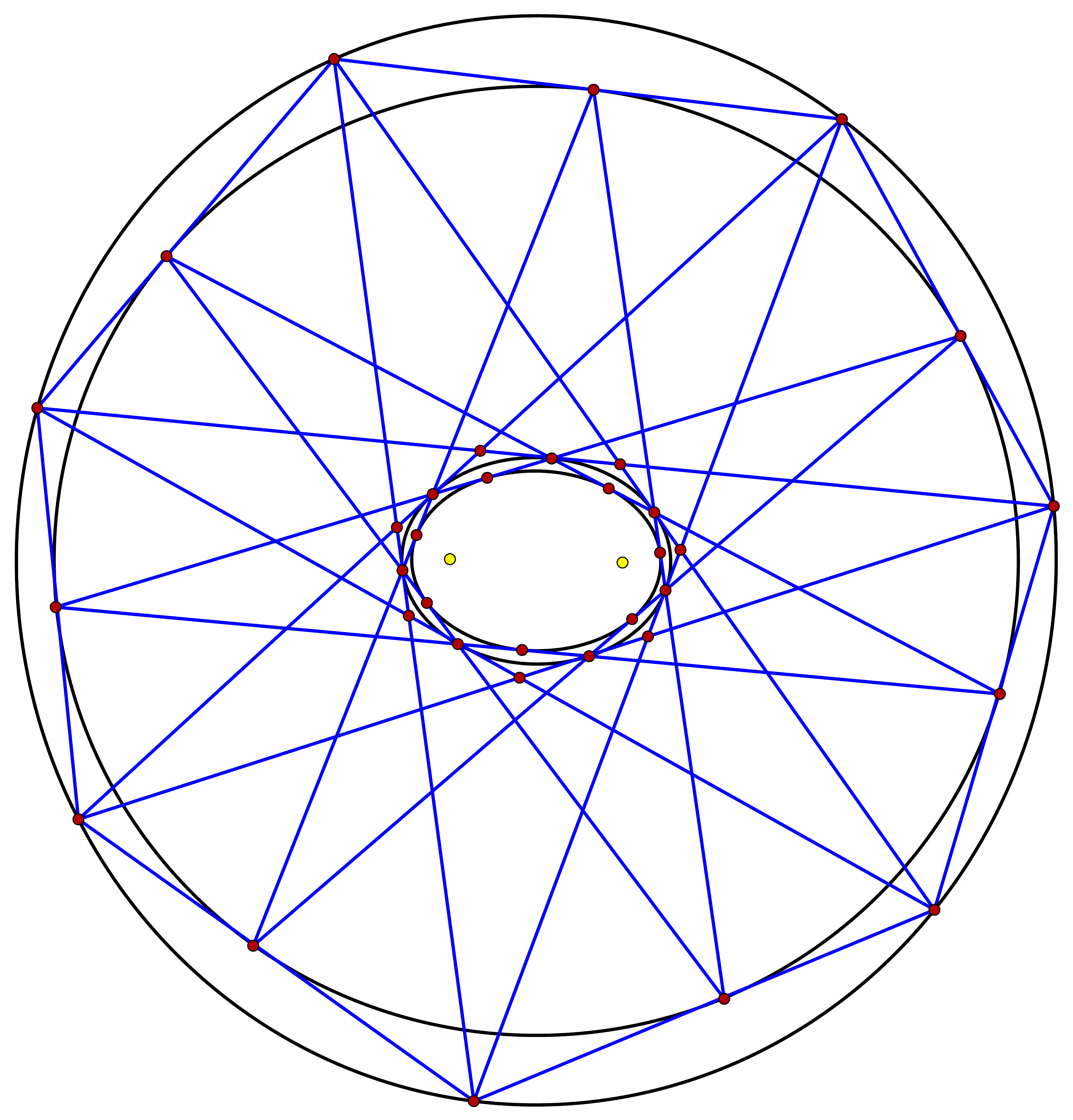}
\caption{Theorem \ref{gen}: $n=7, k=3$.}
\label{hepta}
\end{figure}

\begin{theorem} \label{gen}
The two operations on Poncelet polygons commute: $ID_k(P)=D_kI(P)$, see Figure \ref{hepta}.
\end{theorem}

We shall deduce Theorem \ref{gen} from the Poncelet grid theorem \cite{Sch07} which we now describe. The statements slightly differ for odd and even $n$, and we assume  that $n$ is odd. 

Let $\ell_1,\ldots, \ell_n$ be the lines containing the sides of a Poncelet $n$-gon, enumerated in such a way that their tangency points with the inscribed ellipse are in the cyclic order. The {Poncelet grid} is the collection of $n(n+1)/2$ points  $\ell_i \cap \ell_j$, where  $\ell_i \cap \ell_i$ is, by definition, the tangency point of the line $\ell_i$ with the inscribed ellipse. 

\begin{figure}[hbtp]
\centering
\includegraphics[height=2.2in]{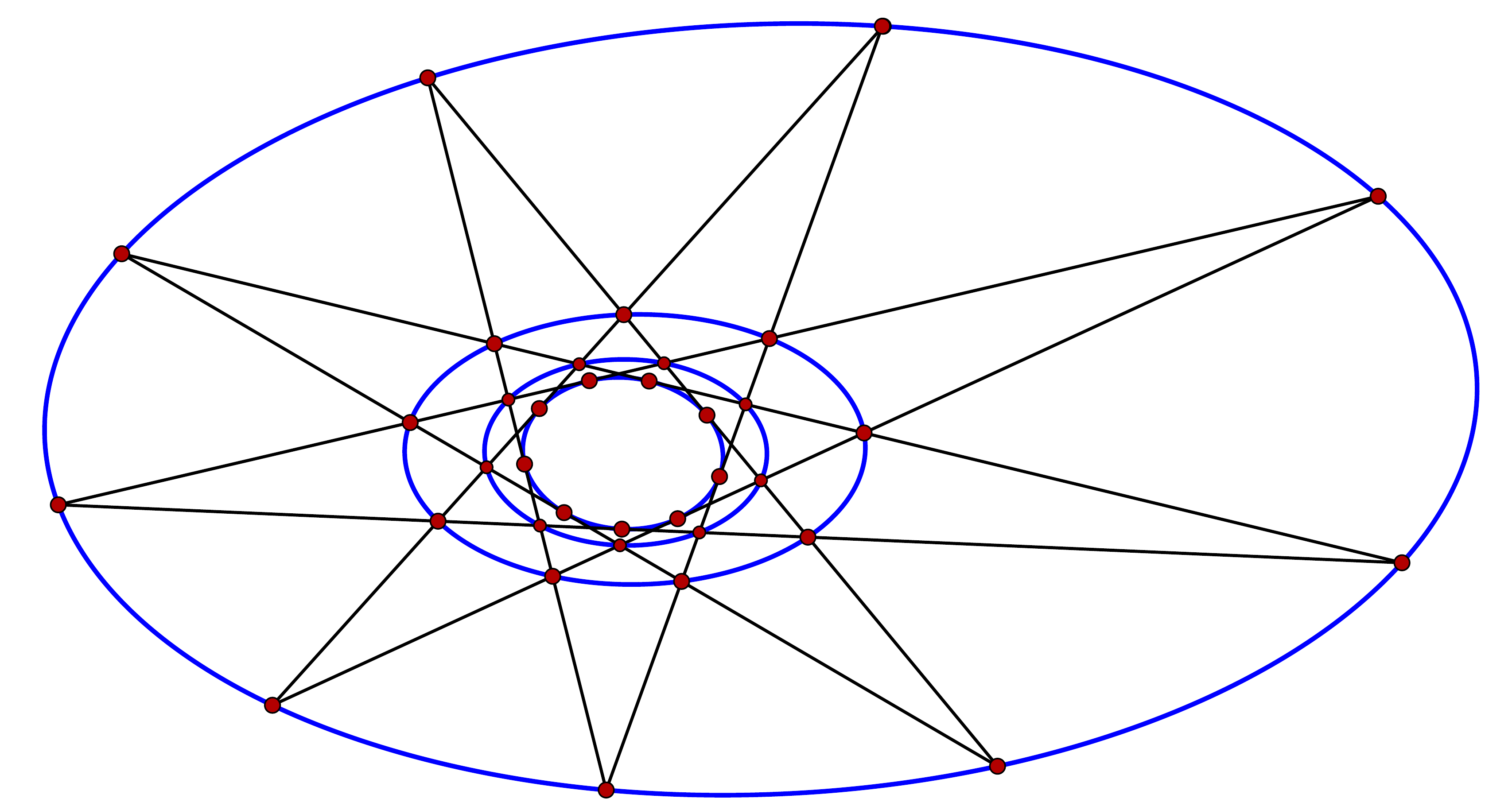}
\caption{Poncelet grid, $n=9$:  the concentric sets $Q_0, Q_2, Q_3$, and $Q_4$  are shown.}
\label{Grid}
\end{figure}

Partition the Poncelet grid in two ways. Define the sets
$$
Q_k = \cup_{i-j=k} \ell_i \cap \ell_j,\quad R_k = \cup_{i+j=k} \ell_i \cap \ell_j,
$$
where the indices are understood mod $n$. One has $(n + 1)/2$ sets $Q_k$ , each containing $n$ points, and $n$ sets $R_k$ , each containing $(n + 1)/2$ points. These sets are called concentric and radial, respectively, see Figure \ref{Grid}. 

\begin{theorem} \label{Pgrid}
(i) The concentric sets lie on nested ellipses, and the radial sets lie on disjoint hyperbolas. \\
(ii) The complexifications of these conics have four common  tangent lines (invisible in Figure \ref{Grid}). \\
(iii) All the concentric sets are projectively equivalent to each other, and so are all the radial sets. 
\end{theorem}

The Poncelet grid theorem is a result in projective geometry. One can apply a projective transformation so that the initial two ellipses that support the Poncelet polygon become confocal (this is the case in Figure \ref{hepta}). Under this assumption, an Euclidean version of the Poncelet grid theorem, proved in \cite{LT}, asserts that all the concentric and the radial sets lie on confocal ellipses and hyperbolas. Furthermore, let 
$$
\frac{x^2}{a^2+\lambda} + \frac{y^2}{b^2+\lambda} =1\quad {\rm and} \quad \frac{x^2}{a^2+\mu} + \frac{y^2}{b^2+\mu} =1
$$
be two ellipses containing concentric sets $Q_i$ and $Q_j$. Then the linear map 
\begin{equation} \label{maps}
A_{\lambda \mu} = \pm {\rm Diag} \left(\sqrt{\frac{a^2+\mu}{a^2+\lambda}}, \sqrt{\frac{b^2+\mu}{b^2+\lambda}}\right),
\end{equation}
that relates these ellipses, takes $Q_i$ and $Q_j$ (the sign depends on the parity of $i-j$). 

Now we prove Theorem \ref{gen}. 
Start with a pair of confocal ellipses $(\Gamma,\gamma)$ and a Poncelet $n$-gon $P$, inscribed into the outer one, $\Gamma$, and circumscribed about the inner one, $\gamma$. In Figure \ref{hepta}, these are the two inner-most ellipses. The tangency points of the sides of $P$ with $\gamma$ form the set $Q_0$, and its vertices form the set $Q_1$. 

Consider the set $Q_k$ (where $k=3$ in Figure \ref{hepta}). These $n$ points lie on a confocal ellipse, say $\delta$ (in Figure \ref{hepta}, this is the third ellipse, counting from within). 

Let $A$ be the linear map of the form (\ref{maps}) that takes $Q_0$ to $Q_1$. This map takes the tangent lines $\ell_i$ to $\gamma$ at the points of the set $Q_0$ to the tangent lines $L_i$ to $\Gamma$ at the points of $Q_1$. Therefore $A$ takes  $\ell_i \cap \ell_{i+k}$ to $L_i \cap L_{i+k},\ i=1,\ldots,n$. 
Call the latter set $S$.  It follows that $S$ lies on an ellipse $A(Q_k) =: \Delta$, and $A(\delta)=\Delta$.

 Let $B$ be the linear map of the form (\ref{maps}) that takes $Q_1$ to $Q_k$. Then the map $BA$ takes $Q_0$ to $Q_k$. The maps of the family (\ref{maps}) commute, hence $BA=AB$. The map $AB$ takes $\Gamma$ to $\Delta$, and it takes the initial Poncelet $n$-gon $P$ to a Poncelet $n$-gon on the ellipses $(\Delta,\delta)$ with the vertex set $S$. This implies the statement of the theorem.
\medskip
 
Note that the ellipse $\Delta$ is not confocal with the confocal ellipses $\gamma, \Gamma$, and $\delta$. In fact, the ellipses $\Delta, \delta$, and $\Gamma$ form a pencil (a 1-parameter family of conics that share four points; these four points are complex and not visible in the figure). Indeed, these three ellipses, along with the heptagons inscribed in $\Delta$ and circumscribed about $\delta$ and $\Gamma$, form a  Poncelet grid projectively dual to the original one.

We finish with some comments placing our result in context.  
 
\paragraph{Remarks}
1)  The map that we denoted by $D_2$ is the well know pentagram map. Considered as a transformation on the moduli space of projectively equivalent polygons, the pentagram map is a completely integrable system, closely related with the theory of cluster algebras; see \cite{GSTV,Gl,GP,OST1,OST2,Sch92,Sch08,So} for some references. The deeper diagonal maps $D_k,\ k>2$, are also completely integrable.

Incidentally, the Poncelet grid theorem implies that the pentagram map sends a Poncelet polygon to a projectively equivalent one. 

2) As a by-product of the study of the pentagram map, eight new projective configuration theorems were found in \cite{ST1}; see also \cite{Ta}.

3) To quote from Kasner \cite{Ka},
\begin{quote}
Among the corollaries ... the most interesting is this (valid for at least convex pentagons):
{\it The limit point of the sequence of successive inscribed pentagons coincides with the limit point of the sequence of successive diagonal pentagons.}
\end{quote}

In this respect, see \cite{Gl2} for a recent striking result on the limiting point of the pentagram map acting on convex $n$-gons (before factorization by the projective group). 

\begin{figure}[hbtp]
\centering
\includegraphics[height=2.2in]{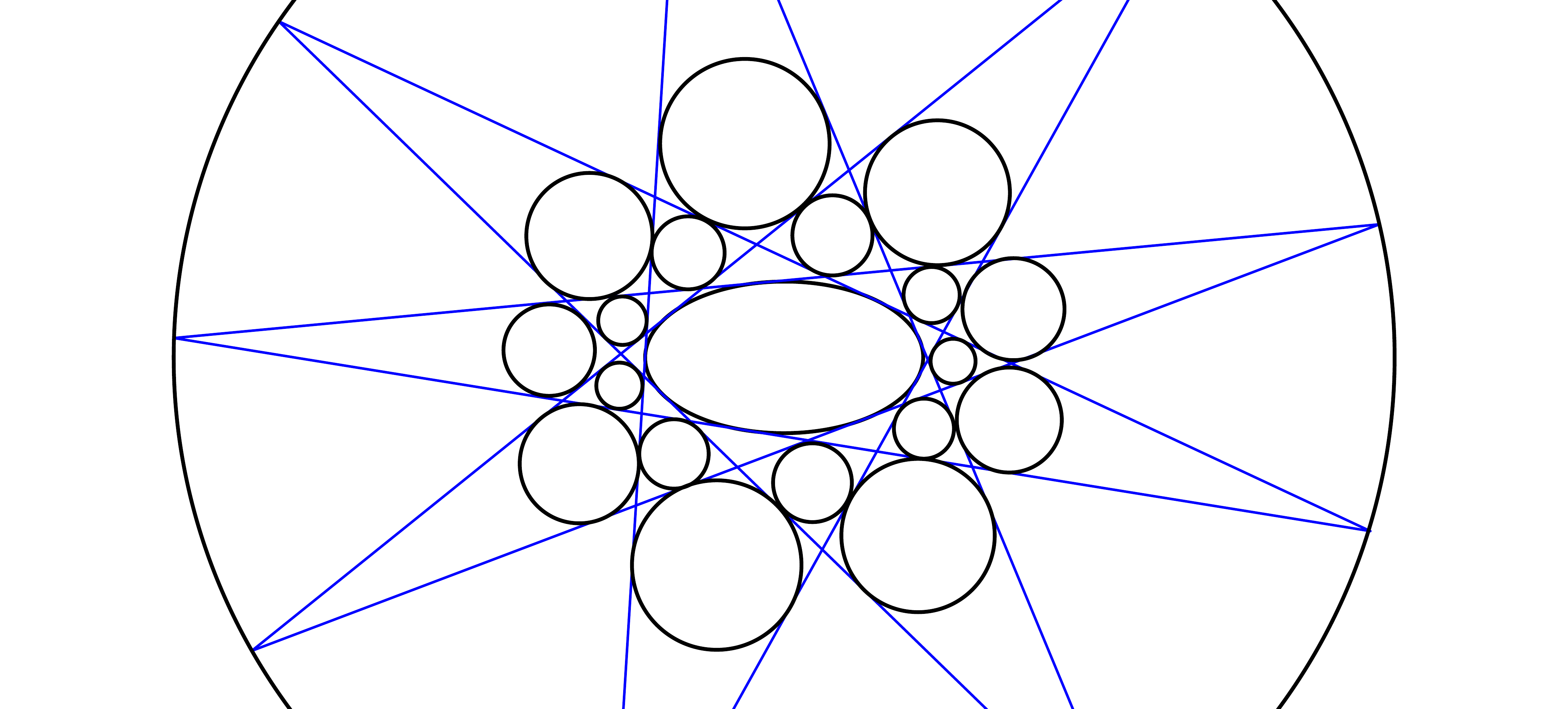}
\caption{Poncelet grid of circles.}
\label{circles}
\end{figure}

4) Poncelet theorem and its ramifications continue to be an active research area. See \cite{AB,ST2} for some recent results. For example, every quadrilateral of the Poncelet based on confocal ellipses admits an inscribed circle \cite{AB}, see Figure \ref{circles}.

5) The name of Edward Kasner (1878--1955) may be lesser known to contemporary mathematicians. He was a prominent geometer who spent his career at Columbia University, see a memoir \cite{Dou} by his student, J. Douglas (one of two winners of the first Fields Medals). We cannot help mentioning that Kasner introduced the terms ``googol" and ``googolplex".\footnote{Just google it!} To quote from the bestselling book \cite{KN} by Kasner and his student Newman,
\begin{quote}
Words of wisdom are spoken by children at least as often as by scientists. The name ``googol" was invented by a child (Dr. Kasner's nine-year-old nephew) who was asked to think up a name for a very big number, namely, 1 with a hundred zeros after it.
\end{quote}

\bigskip
{\bf Acknowledgements}. It is a pleasure to acknowledge stimulating discussions with A. Akopyan, I. Izmestiev, R. Schwartz, and Yu. Suris. 
I was supported by NSF grant DMS-1510055.

\end{document}